\documentclass{amsart}
\usepackage{amssymb}
\usepackage{tikz}
\usepackage{xcolor}
\usepackage{epic,curves,mfpic}
\newtheorem{theo}{Theorem}[section]

\newtheorem{prop}[theo]{Proposition}

\theoremstyle{remark}
\newtheorem{exam}[theo]{Example}

\newcommand{\Z}{\mathbb{Z}}
\newcommand{\R}{\mathbb{R}}

\newcommand{\bv}{{\bf v}}

\newcommand{\bu}{{\bf u}}

\newcommand{\bP}{\mathcal{P}}

\newcommand{\m}{\mu}

\newcounter{fig}
\begin{document}

\title{Proper decompositions of finitely presented groups}
\author{A.N.Bartholomew and M.J.Dunwoody}

\begin{abstract}  This is a report on our long term project to find an algorithm to decide if a finitely presented group has a non-trivial action on a tree.

\end{abstract}
\maketitle

\section{Introduction}

 In his seminal work \cite{[St]} Stallings showed that a finitely generated group with more than one end
 splits over a finite subgroup.   In \cite{[D1]} it was shown that a finitely presented group is accessible. 
 This means that a finitely presented group $G$ has a decomposition as the fundamental group
 of a graph of groups in which vertex groups are at most one ended and edge groups are finite.
 This decomposition provides information about every action of $G$ on a simplicial tree with finite
 edge groups.  Thus,  let $S$ be the Bass-Serre $G$-tree associated with the decomposition described
 and let $T$ be an arbitrary $G$-tree with finite edge stabilizers, then there is a $G$-morphism $\theta :  S \rightarrow T$.
 We say that any action is {\it resolved} by the action on $S$.
 In \cite {[D3]} and \cite {[D4]}  examples are given of inaccessible groups.   These are finitely generated groups - but not 
 finitely presented - for which there is no such $G$-tree $S$.  These groups do have actions on
 a special sort of $\R $-tree (a realization of a protree) but there appears to be no such action
 which resolves all the other actions.
 
 An earlier version of this paper sought to show that a finitely presented group has an action that resolves all actions.  Sadly this is incorrect. The Higman group, discussed below, has two incompatible decompositions  and there is no action on a tree that  that resolves both the trees corresponding to these decompositions.
 
 It is easy to determine if a finitely presented group splits as an HNN-group.  This is the case if and only if the group made abelian is infinite.
 Deciding if a group splits as a free product  with amalgamation is much harder.
 It is known that there
is a group $H$ which has a presentation for which it cannot be decided if the group is non-trivial.
One could use this presentation to construct a presentation for $H*H$.   Clearly it will not be possible
to decide it this decomposition is non-trivial. 

  It seems possible that for a finitely presented group $G$ there is a finite list of decompositions such that if
$G$ has a non-trivial decomposition then a non-trivial decomposition is in this list.  If the group $G$ has a solvable membership algorithm
then it will be possible to decide if a decomposition in the list is non-trivial.
 In an earlier version of this paper we claimed that a list we could construct for $G$ did have the required property.  However our proof was not correct.
 We think that the methods, described here, of determining a list of  different decompositions of a finitely presented group, could yet lead to interesting results.

The  result -and its proof - on the accessibility of finitely presented groups can be seen as a generalization of a result by Kneser (see \cite {[He]}) - and its proof -  that a compact $3$-manifold (without boundary) has a prime decomposition, i.e. it can be expressed
as a connected sum of a finite number of prime factors.  A compact $3$-manifold $M$ is prime
if for every decomposition $M = M_1\sharp M_2$ as a connected  sum, either $M_1$ or $M_2$ is
a $3$-sphere.  Expressed as a result about fundamental groups, it says that the fundamental group
of a compact $3$-manifold is a free product of finitely many factors, which, of course, is true for
any finitely generated group by Grushko's Theorem.

Kneser's result is a basis for the theory of normal surfaces in $3$-manifolds,  due to Haken (see \cite {H}), used to
provide an algorithm to decide if a knot is trivial.    Jaco and Oertel \cite {[JO]} and Jaco and Tollefson \cite {[JT]} used normal surface theory
to develop algorithms for deciding if a compact $3$-maniflold $M$ contains an incompressible surface.    If this is the case then $G = \pi _1M$
splits over a subgroup that is the group of the embedded surface.     The theory of tracks and patterns used in \cite {[DD]} and  \cite {[D1]}
is a generalisation of the theory of normal surfaces.   Instead of using the way a surface intersects the different $3$-simplexes, a pattern is determined by intersections with $2$-simplexes,   A pattern in the $2$-skeleton of a $3$-manifold determines a surface in which the intersection with each $3$-simplex is a finite set of disjoint discs.  This surface 
is called a patterned surface.    The proofs of Jaco and Oertel for normal surfaces will also work for patterned surfaces.   In \cite {[DD]} the theory of patterned surfaces is used to give proofs of the equivariant loop and sphere theorems. It is a natural question to ask if the theory of tracks and patterns  can be used to provide
algorithms for deciding if a finitely presented group splits.   This is because an action of a finitely presented group on a tree is resolved by an action on a tree corresponding to a pattern in the presentation complex,  so that if such a group has a non-trivial action on a tree then there will be a track in the presentation complex giving a non-trivial decomposition.
 In this paper we describe our attempts to answer this question.

A normal surface in a compact $3$-manifold corresponds to a particular solution to a set of matching equations in $\Z ^n$.
These solutions all lie in a polyhedral convex cone in  $\R^n$.    The algorithms referred to above consist of showing that if a there is 
an incompressible surface in $M$, then there is one that corresponds to one in a finite list of points in this cone.    In some cases the list
is just the extreme fundamental solutions, (or vertex solutions)  i.e.  those points that  are the smallest integer valued points in the one dimensionsal faces of the cone.
Thus Jaco and Tollefson show that that there is a face of the cone for which the extreme fundamental solutions give a prime decomposition
of the manifold.    

Using software developed by the first author we have found examples that the results for patterned surfaces in $3$-manifolds cannot be generalised as much as one might hope.
Thus for any finitely presented group $G$, there is a finite $2$-dimensional $2$-complex $X$ with fundamental group $G$.   The tracks in $X$ correspond to  points in a cone $\mathcal P$,  We had been hoping to show that if $G$ has a non-trivial splitting that corresponds to an internal point of the face of $\mathcal P$, then at least one or hopefully all of the vertex solutions of that face will give non-trivial splittings. However this is not always the case.
We give an example in which two trivial vertex solutions have a rational linear combination that gives a non-trivial splitting.

The vertex solutions of $\mathcal P$ and the corresponding decompositions can be computed.   Programmes for doing this are available
on the first author's website.

It still seems likely that if a finitely presented group has a non-trivial decomposition, then there will be a fundamental solution that corresponds to a track giving a non-trival decompostion,  and that there are only finitely many fundamental solutions which lie in a bounded subset of the solution cone.

Here are some questions that remain to be answered. 

 Let $G$ be a finitely presented group, with presentation complex $X$ and corresponding solution cone $\mathcal P$.

1.  If there is a non-trivial homomorphism $G \rightarrow \Z$, then  is there at least one  fundamental solution or even a vertex solution that is non-separating?

2.  If $G$ splits, then is there a non-trivial fundamental solution or even a non-trivial vertex solution?

3.  If $G$ has more than one end, i.e. if $G$ splits over a finite subgroup,  then is there a fundamental solution of even a vertex solution corresponding to a splitting over a finite subgroup?

4.  Do the fuundamental solutions lie in a bounded region of $\mathcal P$?

\section{Tracks and Patterns}

We illustrate the theory by repeated reference to a particular example.

The cell complex for the trefoil group { $G = \langle c, d | c^3 = d^2\rangle $}

  \begin{figure}[htbp]
\centering
\begin{tikzpicture}[scale=1.5]
 \fill[brown]  (1,.3) -- (2, .8) -- (2, 2.2) -- (1,2.7) -- (0,1.5) -- cycle ;
  \draw  (1,.3) -- (2, .8) -- (2, 2.2) -- (1,2.7) -- (0,1.5) -- (1, .3) ;
  \draw  (.55, 2.25) --   (.5, 2.1) -- (.64, 2.18) ;
    \draw  (.45,1.05) --   (.5, .9) -- (.36, .96) ;
  \draw  (1.65, .56) --   (1.5, .55) -- (1.61, .68) ;
  \draw  (1.35, 2.46) --   (1.5, 2.45) -- (1.39, 2.58) ;
  \draw  (2, 1.7) node {$\vee $} ;
  \draw  (2.2, 1.7) node {$c$} ;
\draw  (1.7, .45) node {$c$} ;
\draw  (.4, .8) node {$d$} ;
\draw  (.4, 2.2) node {$d$} ;
\draw  (1.6, 2.6) node {$c$} ;
\draw  (4.25, 2.5) node {$d$} ;
\draw  (4.25, .5) node {$c$} ;
\draw  (4.25, 2.25) node {$>$} ;
\draw  (4.25, .75) node {$>$} ;

\draw (4, 1) -- (4.5, 2) ;
\draw (4.5, 1) -- (4, 2) ;
\draw (4.5, 2) arc (0: 180 : .25 cm) ;
\draw (4, 1) arc (180: 360 : .25 cm) ;

\end{tikzpicture}
\end{figure}
Attach the $5$-sided disc to the figure eight as specified by the letters and arrows.
The space $X$ has { $\pi _1(X) = G$.}

A group presentation can be changed so that every relation has length at most three, giving a presentation complex with $2$-cells having  at most
$3$ edges.

 \begin{figure}[htbp]
\centering
\begin{tikzpicture}[scale=1.5]
 \fill[brown]  (1,.3) -- (2, .8) -- (2, 2.2) -- (1,2.7) -- (0,1.5) -- cycle ;

  \draw  (1,.3) -- (2, .8) -- (2, 2.2) -- (1,2.7) -- (0,1.5) -- (1, .3) ;
  \draw  (.55, 2.25) --   (.5, 2.1) -- (.64, 2.18) ;
    \draw  (.45,1.05) --   (.5, .9) -- (.36, .96) ;
  \draw  (1.65, .56) --   (1.5, .55) -- (1.61, .68) ;
  \draw  (1.35, 2.46) --   (1.5, 2.45) -- (1.39, 2.58) ;
  \draw  (2, 1.7) node {$\vee $} ;
  \draw  (2.2, 1.7) node {$c$} ;
\draw  (1.7, .45) node {$c$} ;
\draw  (.4, .8) node {$d$} ;
\draw  (.4, 2.2) node {$d$} ;
\draw  (1.6, 2.6) node {$c$} ;
\draw  (4.25, 2.5) node {$d$} ;
\draw  (4.25, .5) node {$c$} ;
\draw  (4.25, 2.25) node {$>$} ;
\draw  (4.25, .75) node {$>$} ;
\draw  (.8, 1.5) node {$e$} ;
\draw  (3.3, 1.5) node {$e$} ;
\draw  (1.6, 1.8) node {$f$} ;
\draw  (5.2, 1.5) node {$f$} ;

\draw (4, 1) -- (4.5, 2) ;
\draw (4.5, 1) -- (4, 2) ;
\draw (4.5, 2) arc (0: 180 : .25 cm) ;
\draw (4, 1) arc (180: 360 : .25 cm) ;
\draw [red] (1, .3) -- (1, 2.7) ;
   \draw [red] (2, .8) -- (1, 2.7) ;
\draw [red] (3.75, 1.25) -- (4.75, 1.75) ;
\draw [red]  (3.75, 1.75) -- (4.75, 1.25) ;
\draw [red] (3.75, 1.75) arc (90: 270 : .25 cm) ;
\draw [red] (4.75, 1.25) arc (-90: 90: .25 cm) ;
 \draw  (3.5, 1.5) node {$\vee $} ;
 \draw  (5, 1.5) node {$\vee $} ;
 \draw  (1, 1.5) node {$\vee $} ;
  \draw  (1.5, 1.9) --   (1.5, 1.75) -- (1.38, 1.8) ;

\end{tikzpicture}

\end{figure}

Thus {$G  =\langle c, d | c^3 = d^2 \rangle =  \langle c, d, e, f | d^2 = e, e =fc, f =c^2\rangle $.}

The cell complex $X$ consists of three $3$-sided $2$-cells  attached to a $4$-leaved rose.

\medskip

Let $X$ be a cell complex in which each $2$-cell is $3$-sided.

A {\it pattern} is a subset of $X$ which intersects each $2$-cell in a finite number of disjoint lines each of which
intersects the boundary of the $2$-cell in its two end  points which lie in distinct edges.

A {\it track} is a connected pattern.

If $X$ has $m$ $2$-cells then a pattern is specified (up to an obvious equivalence) by
a $3m$-vector  in which there are three coefficients for each $2$-cell which record the number of   lines joining 
the two edges at each corner.   Thus for the $2$-cell

  \begin{figure}[htbp]
\centering
\begin{tikzpicture}[scale=1]

  \draw  (0,.0) -- (4, 0) -- (2, 4) -- (0,0)  ;
{

\draw [red] (1.5, 3) -- (2.25, 3.5) ;
\draw [red]  (1.25, 2.5) -- (2.75, 2.5) ;
\draw [red]  (.5, 0) -- (.4, .8) ;
\draw [red]  (.7, 0) --(.8, 1.6);
\draw [red]  (1.1, 0) -- (1, 2) ;
\draw [red] (3, 0) -- (3.5, 1) ;
\draw [red]  (2.5, 0) -- (3, 2) ;

}
\end{tikzpicture}

\flushleft the coefficients $2, 2, 3 $ record the intersection of the pattern with the $2$-cell.
\end{figure}

\def\x{{\bf x}}
\def\y{{\bf y}}
\def\t{{\bf t}}
\def\p{{\bf p}}
\def\q{{\bf q}}
\def\u{{\bf u}}
\def\m{{\bf m}}
\def\n{{\bf n}}
\def\1{{\bf 1}}

\def\a{{\bf a}}

\def\v{{\bf v}}
\def\s{{\bf s}}
\def\w{{\bf w}}

For the complex $X$ for the trefoil group $G$
a pattern is specified by a $9$-vector, where the $i$-th coefficient corresponds to the number of lines 
crossing the $i$-th corner  labelled $i$ in red in the diagram below.
  \begin{figure}[htbp]
\centering
\begin{tikzpicture}[scale=1.5]

  \draw  (1,.3) -- (2, .8) -- (2, 2.2) -- (1,2.7) -- (0,1.5) -- (1, .3) ;
  \draw  (.55, 2.25) --   (.5, 2.1) -- (.64, 2.18) ;
    \draw  (.45,1.05) --   (.5, .9) -- (.36, .96) ;
  \draw  (1.65, .56) --   (1.5, .55) -- (1.61, .68) ;
  \draw  (1.35, 2.46) --   (1.5, 2.45) -- (1.39, 2.58) ;
  \draw  (2, 1.7) node {$\vee $} ;

 \draw [dashed] (1, .3) -- (1, 2.7) ;
   \draw [dashed] (2, .8) -- (1, 2.7) ;
   {
\draw [red] (1.2, 2.5) node {$_1$} ;
\draw [red] (1.8, 2.2) node {$_2$} ;
\draw [red] (1.9, 1.2) node {$_3$} ;
\draw [red] (1.1, 2.3) node {$_4$} ;
\draw  [red] (1.8, .9) node {$_5$} ;
\draw  [red] (1.1,.5) node {$_6$} ;
\draw [red]  (.9, 2.4) node {$_7$} ;
\draw [red]  (.9, .55) node {$_8$} ;
\draw  [red] (.2, 1.5) node {$_9$} ;}
\draw  (1.6, 2.6) node {$c$} ;
\draw  (2.2, 1.7) node {$c$} ;
\draw  (1.7, .45) node {$c$} ;
\draw  (.4, .8) node {$d$} ;
\draw  (.4, 2.2) node {$d$} ;
\draw  (1.5, 1.75) node {$f$} ;
\draw  (.9,1.5) node {$e$} ;
\end{tikzpicture}

\end{figure}
\bigskip
  In the trefoil complex  a  vector of non-negative integers  { $x = (x_1, x_2, \dots , x_9)$}  is a pattern in if it satisfies  the matching equations
$$x_1 +x_2 = x_2 +x_3 = x_5+ x_6$$ (number of intersection points with edge $c$)

$$x_1 + x_3 = x_ 4 + x_5$$(number of intersection points with edge $f$)

$$x_4 + x_6 = x_7 +x_8$$(number of intersection points with edge $e$)

$$x_ 7 +x_9 = x _8 + x_9$$ (number of intersection points with edge $d$)

In general a  $3m$-vector corresponds to a pattern, if and only if
\begin{itemize} 
\item[{(i)}] Each entry is a non-negative integer.
\item[{(ii)}] It is a solution vector to a finite set of linear equations called the {\it matching equations}, where if an edge $e$ lies in
$k$ $2$-simplexes, then there are $k-1$ matching equations corresponding to the intersection of the pattern with $e$.
\end{itemize}

In general a pattern $P$ in a $2$-complex $X$ will lift to a pattern $\tilde P$  in $\tilde X$.  Each track component of $\tilde P$ will separate
and there is a $G$-tree $T_P$ in which the edges correspond to the track components of $\tilde P$ (see \cite {[DD]}, Chapter VI or \cite {[D1]} 
for details).   If $P$ consists of a single track then $T_P$ will be the Bass-Serre tree for a decomposition of $G$ as a free product with amalgamation,
if the track is separating, and as an HNN-group if it is untwisted and  non-separating.   An {\it untwisted } track $t$ is one which has 
a neighbourhood that is homeomorphic to $t \times I$ where $I$ is a closed interval.

  \begin{figure}[htbp]
\centering
\begin{tikzpicture}[scale=1.5]

  \draw  (1,.3) -- (2, .8) -- (2, 2.2) -- (1,2.7) -- (0,1.5) -- (1, .3) ;
  \draw  (.55, 2.25) --   (.5, 2.1) -- (.64, 2.18) ;
    \draw  (.45,1.05) --   (.5, .9) -- (.36, .96) ;
  \draw  (1.65, .56) --   (1.5, .55) -- (1.61, .68) ;
  \draw  (1.35, 2.46) --   (1.5, 2.45) -- (1.39, 2.58) ;
  \draw  (2, 1.7) node {$\vee $} ;

 \draw [dashed] (1, .3) -- (1, 2.7) ;
   \draw [dashed] (2, .8) -- (1, 2.7) ;
   {
\draw  (1.2, 2.5) node {$_1$} ;
\draw  (1.8, 2.2) node {$_1$} ;
\draw  (1.9, 1.2) node {$_1$} ;
\draw  (1.1, 2.3) node {$_0$} ;
\draw  (1.8, .9) node {$_2$} ;
\draw  (1.1,.5) node {$_0$} ;
\draw  (.9, 2.4) node {$_0$} ;
\draw  (.9, .55) node {$_0$} ;
\draw  (.2, 1.5) node {$_0$} ;}
\draw  (1.6, 2.6) node {$c$} ;
\draw  (2.2, 1.7) node {$c$} ;
\draw  (1.7, .45) node {$c$} ;
\draw  (.4, .8) node {$d$} ;
\draw  (.4, 2.2) node {$d$} ;

\end{tikzpicture}

\end{figure}

  \begin{figure}[htbp]
\centering
\begin{tikzpicture}[scale=1.5]

  \draw  (1,.3) -- (2, .8) -- (2, 2.2) -- (1,2.7) -- (0,1.5) -- (1, .3) ;
  \fill[green] (1.2, .4) -- (1.2, 2.6) -- (1.8, 2.3) -- (2, 2) --(2,1) -- (1.8, .7) -- cycle ;

  \draw  (.55, 2.25) --   (.5, 2.1) -- (.64, 2.18) ;
    \draw  (.45,1.05) --   (.5, .9) -- (.36, .96) ;
  \draw  (1.65, .56) --   (1.5, .55) -- (1.61, .68) ;
  \draw  (1.35, 2.46) --   (1.5, 2.45) -- (1.39, 2.58) ;

  \draw  (2, 1.7) node {$\vee $} ;

 \draw [dashed] (1, .3) -- (1, 2.7) ;
   \draw [dashed] (2, .8) -- (1, 2.7) ;
   
\draw [blue] (1.2, .4) -- (1.2, 2.6) ;
\draw [blue] (1.8, 2.3) -- (2, 2) ;
\draw [blue] (1.8, .7) -- (2, 1) ;

\draw  (1.6, 2.6) node {$c$} ;
\draw  (2.2, 1.7) node {$c$} ;
\draw  (1.7, .45) node {$c$} ;
\draw  (.4, .8) node {$d$} ;
\draw  (.4, 2.2) node {$d$} ;

\end{tikzpicture}

\end{figure}

In the trefoil complex $X$ an example of a pattern   is as follows.     The $9$-vector 

{ $ t  _1 = (1, 1, 1, 0, 2, 0, 0, 0, 0 )$}

corresponds to the pattern shown above.  Thus there is one line crossing each of the corners labelled $1, 2$ and $3$
and $2$ lines crossing the corner labelled $5$.

This pattern  is in fact a separating track and corresponds to the decomposition of $G$.

$$G = \langle d \rangle *_{\langle d^2 = c^3\rangle } \langle c\rangle .$$

The track separates into two regions one of which is coloured green.

A separating track is always untwisted.   It $t$ is twisted, then $2t$ is separating and hence untwisted.

 The   track $t$ shown below in blue  is  { twisted}   so the pattern $2t$ is also a track. The separating track $2t$ gives the trivial decomposition 
 $G = G*_HH$ where $H$ has index two in $G$

   \begin{figure}[htbp]
\centering
\begin{tikzpicture}[scale=2]
  \draw  (1,.3) -- (2, .8) -- (2, 2.2) -- (1,2.7) -- (0,1.5) -- (1, .3) ;
  \draw  (.55, 2.25) --   (.5, 2.1) -- (.64, 2.18) ;
    \draw  (.45,1.05) --   (.5, .9) -- (.36, .96) ;
  \draw  (1.65, .56) --   (1.5, .55) -- (1.61, .68) ;
  \draw  (1.35, 2.46) --   (1.5, 2.45) -- (1.39, 2.58) ;
  \draw  (2, 1.7) node {$\vee $} ;

 \draw [dashed] (1, .3) -- (1, 2.7) ;
   \draw [dashed] (2, .8) -- (1, 2.7) ;
   {
  \draw [blue] (.6, .8) -- (1.6, .6) ;
  \draw[blue]  (.6,2.2) -- (1.8, .7) ;
  \draw [blue] (2, 2) -- (1.8, 2.3) ;
\draw [blue] (2, 1.8) -- (1.6, 2.4) ;

}
\draw  (1.6, 2.6) node {$c$} ;
\draw  (2.2, 1.7) node {$c$} ;
\draw  (1.7, .45) node {$c$} ;
\draw  (.4, .8) node {$d$} ;
\draw  (.4, 2.2) node {$d$} ;

\end{tikzpicture}

\end{figure}

The track shown in red is non-separating and untwisted, and gives a decomposition of $G$ as an HNN-group.

  \begin{figure}[htbp]
\centering
\begin{tikzpicture}[scale=1.5]
  \draw  (1,.3) -- (2, .8) -- (2, 2.2) -- (1,2.7) -- (0 ,1.5) -- (1, .3) ;
{  \draw  (.55, 2.25) --   (.5, 2.1) -- (.64, 2.18) ;
    \draw  (.45,1.05) --   (.5, .9) -- (.36, .96) ;
\draw  (1.65, .56) --   (1.5, .55) -- (1.61, .68) ;
  \draw  (1.35, 2.46) --   (1.5, 2.45) -- (1.39, 2.58) ;
  \draw  (2, 1.7) node {$\vee $} ;}

{ \draw [red]  (1.2, .4) -- (.85, .5) ;
 \draw [red] (1.8, .7) -- (.66, .7) ;

 \draw [red] (2, 1.2) -- (.57, .8) ;
 \draw [red] (2,1.8) -- (0.18 , 1.7 ) ;
  
   \draw [red] (.32, 1.9) --  (1.8,2.3) ;
 \draw [red] (.5, 2.1) -- (1.2, 2.6) ;}

 \draw [dashed] (1, .3) -- (1, 2.7) ;
   \draw [dashed] (2, .8) -- (1, 2.7) ;

{
\draw  (1.6, 2.6) node {$c$} ;
\draw  (2.2, 1.7) node {$c$} ;
\draw  (1.7, .45) node {$c$} ;
\draw  (.4, .8) node {$d$} ;
\draw  (.4, 2.2) node {$d$} ;}

\end{tikzpicture}

\end{figure}

Such a track is always associated with a homomorphism $G \rightarrow \Z$.  In this case
$c \mapsto 2,d\mapsto 3$.

If $X$ has $n \  2$-simplexes and $m \ 1$-simplexes (edges) then $X_1$ has $3n \ 2$-cells and 
$3n + m\ 1$-cells.  A {\it marking} of $X_1$ is a solution to the matching equations.  A marking will be  any point of   a compact, convex linear cell in $\R ^{3n +m}$
called the projective solution space $\bP$.
This theory is a generalization of the theory of normal surfaces or patterned surfaces
in $3$-manifolds (see \cite{[JO]},\cite{[JT]}  and  \cite{[DD]}, Chapter VI).  The {\it extreme} or {\it vertex}  solutions are the ones corresponding to
vertices of the projective solution space. Jaco-Oertel \cite{[JO]} and Jaco-Tollefson \cite{[JT]} have shown that
vertex solutions carry important information about normal surfaces in a $3$-manifold. Thus 
in \cite{[JT]} it is shown that there is a face of $\mathcal P$ for which the vertex solutions give a set of 
$2$-spheres giving a complete factorization of a closed $3$-manifold.
A solution is a vertex solution $\bv $ if it has integer coefficients and integer multiples of $\bv$ are the only solutions to  $n\bv = \bv _1 +\bv_2 $, where $n$ is a positive integer and $\bv_1, \bv _2$
are non-zero vectors in $\mathcal P$  with non-negative integer coefficients. 
The first  author, in his D.Phil. Thesis  \cite{[B]}  investigated the solution
space for a group presentation on a computer.   It was hoped to show that
at least one  vertex solution gives a non-trivial decomposition if the group has such a  decomposition.
We are  still unable to show that this is the case.   It is the case in all the examples we have investigated, but we have counterexamples to stronger results we thought might be true.

 Two patterns  are $\it equivalent$ if they have the same number of intersections with each
 edge, so that they determine the same vector $\bu $.   Two tracks $t_1, t_2$  are {\it compatible }
 if there is a pattern with two components which are equivalent to $t_1$ and $t_2$.
 A track is a {\it fundamental solution} if it cannot be written as a sum of more than one track.  Clearly vertex solutions are fundamental solutions.

  Each separating  track gives a decomposition
of $G$ as a free product with amalgamation (possibly trivial).  Each non-separating track gives a decomposition of $G$ as an HNN-group.

For the trefoil example  the software developed by the first author gives the following output.

{\tt\small\indent
G = < c, d | ccc = dd>

There are five vertex solutions.

Vertex solutions (extreme fundamental tracks), n=9 s=5

  1.   1  1  1  0  2  0  0  0  0
  
  2.   0  0  0  0  0  0  0  0  1
  
  3.   0  2  0  0  0  2  1  1  0
  
  4.   2  0  2  2  2  0  1  1  0
  
  5.   2  0  2  4  0  2  3  3  0
}
\medskip

The first vertex track is the one illustrated above as $t_1$.  

The second vertex track $t_2$ is twisted.  It has a neighbourhood that is a M\" obius Band.  The programme gives the decomposition corresponding to $2t_2$,
The boundary of the M\" obius Band, which is a separating track giving a non-trivial decomposition.

The third vertex solution $t_3$ is also twisted, and is illustrated above as the blue track and $2t_3$ gives a trivial decomposition.

The fourth vertex solution is similar to the third.

The fifth vertex solution $t_5$ is the one illustrated above as the red track.

The track $t_5$ is  non-separating and untwisted.

A different presentation of the trefoil group shows interesting behaviour of tracks.

We first state an easily proved   result about twisted tracks.
A track $t$ is untwisted  if $2t$ is a pattern consisting of two copies of $t$.  If $t$ is twisted, then $2t$ is a separating untwisted track, so that $4t$ is a pattern consisting of two copies of $2t$.
\begin {prop}\label {twist}  Let $t$ be a twisted track. There are two possibilities for the decomposition of $G$ associated with $2t$.
\begin {itemize} 
\item [(i)] The decomposition is trivial. One  vertex group is $G$.  The other vertex group and the edge group are both a subgroup of index $2$ in $G$.
\item [(ii)] The decomposition is non-trivial and the edge group has index $2$ in one of the vertex groups.
\end {itemize}
\end {prop}
  \begin{figure}[htbp]
\centering
\begin{tikzpicture}[scale=1.5]

  \draw  (1,.3) -- (2, .8) -- (2, 2.2) -- (1,2.7) -- (0,2.2) -- (0, .8) --(1, .3) ;
  \draw  (1.65, .56) --   (1.5, .55) -- (1.61, .68) ;
  \draw  (1.35, 2.46) --   (1.5, 2.45) -- (1.39, 2.58) ;
  \draw  (2, 1.7) node {$\vee $} ;
 \draw  (0, 1.7) node {$\vee $} ;
 \draw  (.35, .56) --   (.5, .55) -- (.39, .68) ;
  \draw  (.65, 2.46) --   (.5, 2.45) -- (.61, 2.58) ;

\draw [dashed] (0,.8) -- (1, 2.7) ;

 \draw [dashed] (1, .3) -- (1, 2.7) ;
   \draw [dashed] (2, .8) -- (1, 2.7) ;
   
\draw [red] (.8, 2.5) node {$_{10}$} ;
\draw [red] (.2, 2.2) node {$_{12}$} ;
\draw [red] (.1, 1.3) node {$_{11}$} ;
\draw [red] (1.1, 2.3) node {$_4$} ; 
\draw  [red] (1.8, .9) node {$_5$} ;
\draw  [red] (1.1,.5) node {$_6$} ;
\draw [red]  (.9, 2.3) node {$_7$} ;
\draw [red]  (.9, .55) node {$_8$} ;
\draw  [red] (.2, .9) node {$_9$} ;
\draw  (1.6, 2.6) node {$a$} ;
\draw  (2.2, 1.7) node {$b$} ;
\draw  (1.7, .45) node {$a$} ;
\draw  (.4, 2.6) node {$b$} ;
\draw  (-.2, 1.7) node {$a$} ;
\draw  (.3, .45) node {$b$} ;
\draw [red] (1.2, 2.5) node {$_1$} ;
\draw [red] (1.8, 2.2) node {$_2$} ;
\draw [red] (1.9, 1.2) node {$_3$} ;

s1,s10
\end{tikzpicture}

\end{figure}

An alternative presentation $B = <a, b | aba = bab >$ for the trefoil group provides a number of examples in which  what one might have hoped to be correct turns out to be not the case.
A pattern for this presentation will be determined by a $(12)$--tuple. Where the entries in the $(12$)-tuple are given by the number of lines crossing the corners as in the diagram above.
Note that there is an automorphism $\alpha $ of $B$ that transposes $a$ and $b$, and $\alpha $ induces an autoomorphism of the cell complex and also of the solution space 
$\bP$, which permutes the entries in each $(12)$-tuple by the permutation $(1,10)(2,12)(3,11)(4,7)(5,9)(6,8)$.
For this presentation, there are $15$ vertex solutions $s1, s2, \dots , s15$.  The automorphism $\alpha $ induces the permutation 
$$ (s1, s6)(s2,s14)(s4,s9)(s5,s8)(s7,s13)(s12,s15)$$
The vertex solutions $s3, s10,s11$ are all fixed by $\alpha$.

We have

$s1 = (1, 0, 1, 1, 1, 0, 0, 1, 0, 0 ,0 ,1) $ is a twisted track as in Proposition \ref {twist} (ii) so that $2s1$ gives a separating track giving a non-trivial decomposition in which one factor is generated by $ba$. and the other by $aba$.
The vertex tracks $s3, s4, s7, s9, s11,  s13$  are twisted tracks as in Proposition \ref {twist} (i).  Thus
 $2s3 = (0 ,2 ,0, 0 ,0 , 2, 0, 2 ,0, 0 ,0 ,2) $,  $2s4 = (4,2,0,4,0,6,8,2,0,2,6,0)$, $2s7 = (4,2,0,0,4,2,2,0,2,0, 4,2) $, 

$2s11 = ((2,0,2,2,2,0,2,0,2,2,2,0)$. and $2s12 = (6,0,2,2,6,0,2,0,2,0,4,2)$ 
have trivial decompositions in which one vertex is $G$ and the other has index two in $G$.

The vertex tracks $s2 = (3,1,1,0,4,0,0,0,2,0,2,2), s5 = (2,0,0,1,1,1,2,0,0,0,2,0)$ and their images $s14 = (0,2,2,0,2,0,0,0,4,3,1,1)$ and $s8= (0,0,2,2,0,0,1,1,1,2,0,0) $
are all untwisted tracks giving trivial decompositions.  
Finally 

$s10 =(1,0,1,2,0,1,2,1,0,1,1,0)$ is untwisted and non-separating, and so it gives a decomposition of $G$ as an HNN- extension.   The vertex group is the kernel of the 
homomorphism to $\Z$ in which both $a$ and $b$ are mapped to $1$.  Note that $s10$ is the only vertex solution that is untwisted and non-trivial.
Note that there is no vertex track that is untwisted and separating and corresponds to the non-trivial decomposition.

We have the interesting relation
$$s2 + s14 = 3(1,1,1,0,2,0,0,0,2,1,1,1) $$
where $f = (1,1,1,0,2,0,0,0,2,1,1,1)$ is a track that is separating and untwisted.  It is a non-trivial fundamental solution but not a vertex solution.

The track $f$ is compatible with both $x2$ and $s14$, even though $s2$ and $s14$ are not compatible.  This means that we have the following relations for positive integers
$m, n$ where $n \geq m$

$$ms2 +ns14= 3mf + (n-m)s2, $$
 $$ns2 +ms14 = 3mf + (n-m)s14.$$

We had been hoping that if a group $G$ had a non-trivial splitting then it would show up as a vertex solution.  This is not the case with this presentation of the trefoil group.
Thus $2s1$ and $f$ give the splitting as a free product with amalgamation, but no vertex solution does give this splitting.
The tracks $2s1$ and $f$ are compatible.
The tracks $s5$ and $s8 = \alpha s5$ have contrasting behaviour to $s2$ and $s14 =\alpha s2$.  In this case if $m, n$ are coprime positive integers, then
$ms5 +ns8$ is a track giving a trivial decomposition, or at least looking at a lot of cases suggests that this is the case.

\section{Computing decompositions}
A programme is available on the first author's website that calculates the extreme fundamental solutions (or vertex solutions) for 
  the presentation complex of a finitely 
presented group $G$.  The programme then calculates the decomposition 
corresponding to each such track and identifies those that are clearly 
trivial.  The remaining decompositions are left for manual inspection.  Usually there are more trivial decompositions.

See  http://www.layer8.co.uk/maths/tracks.htm

We present some  output for the Higman group.

\begin {exam}\label {Hig}
Let $H = \langle  a, b, c, d | aba^{-1} = b^2, bcb^{-1} = c^2, cdc^{-1} = d^2, dad^{-1} = a^2\rangle $.

This group was investigated by Higman \cite {[H]}.
He showed that it was infinite and had no non-trivial finite homomorphic images.
His proof that it was non-trivial involved showing that it had a decomposition as a free product with amalgamation
$$H = \langle a, b, c \rangle * _{\langle a, c \rangle}  \langle a, d, c\rangle .$$

Also $\langle a, b, c\rangle $ is the free product with amalgamation 
 
 $$\langle a, b, c\rangle = \langle a, b\rangle *_{\langle b\rangle} \langle b,c\rangle , $$
 where both $\langle a, b\rangle $ and $\langle b, c \rangle$ are isomorphic to the Baumslag-Solitar group $BS (1,2)$.

 \end{exam}

For this group presentation, there are 1429 vertex solutions.  All but 4 of these solutions give trivial decompositions.   The ones giving non-trivial decompositions are numbered 
1, 2, 7 and 739.    The output for the first two tracks is as follows.
\medskip

{\tt\small\indent\obeylines

Group presentation:

  a  b  c  d : ab-a-b-b bc-b-c-c cd-c-d-d da-d-a-a 

Jobname: higman4

\medskip

extreme fundamental track 1

The separating track

(2, 0, 0, 0, 2, 0, 0, 0, 0, 0, 0, 0, 0, 0, 0, 0, 0, 0, 0, 0, 0, 0, 0, 0, 0, 0, 0, 0, 0, 2, 2, 0, 0, 1, 1, 1)

Edge stabilizer generators.
a-b-a
b
aad-a
d

First vertex stabilizer generators.
ab-a
aad-a
b
d
c

Second vertex stabilizer generators.
b
a
d

\medskip

extreme fundamental track 2

The separating track

(0, 0, 0, 0, 0, 0, 0, 0, 0, 0, 0, 0, 0, 0, 0, 0, 0, 0, 0, 0, 2, 2, 0, 0, 1, 1, 1, 2, 0, 0, 0, 2, 0, 0, 0, 0)

Edge stabilizer generators.
ddc-d
c
d-a-d
a

First vertex stabilizer generators.
d-c-d-d
d-a-d
c
a
b

Second vertex stabilizer generators.
d
c
a

}

\medskip

In fact it seems these are the only tracks giving non-trivial decompositions.  Taking linear combinations of an incompatible pair of these non-trivial decompostions only appears to produce trivial decompositions.

\medskip

{\tt\small\indent\obeylines

Adding patterns

Tracks from .eft file:

Track 1.    2   0   0   0   2   0   0   0   0   0   0   0   0   0   0   0   0   0   0   0   0   0   0   0   0   0   0   0   0   2   2   0   0   1   1   1

Track 739.    0   0   2   2   0   0   1   1   1   2   0   0   0   2   0   0   0   0   0   0   0   0   0   0   0   0   0   0   0   0   0   0   0   0   0   0

Sum of patterns is the pattern:
2 0 2 2 2 0 1 1 1 2 0 0 0 2 0 0 0 0 0 0 0 0 0 0 0 0 0 0 0 2 2 0 0 1 1 1 
which is a track

\bigskip

Decomposing a given track.

The separating track

(2, 0, 2, 2, 2, 0, 1, 1, 1, 2, 0, 0, 0, 2, 0, 0, 0, 0, 0, 0, 0, 0, 0, 0, 0, 0, 0, 0, 0, 2, 2, 0, 0, 1, 1, 1)

Gives a trivial decomposition.

\medskip

Edge stabilizer generators.
a-b
aba-b-a
ab-a-bcba-b-a
c
a-d-a-a
d

\medskip

First vertex stabilizer generators.
aad-a
aba-b-a
ab-a-b-cba-b-a
a-b
d
c

\medskip

Second vertex stabilizer.
G
}

\end{document}